\documentclass[12pt]{amsart}
\usepackage{fullpage,url,amssymb,enumerate,colonequals}
\usepackage{mathrsfs} 

\usepackage[OT2,T1]{fontenc}
\DeclareSymbolFont{cyrletters}{OT2}{wncyr}{m}{n}
\DeclareMathSymbol{\Sha}{\mathalpha}{cyrletters}{"58}

\usepackage{color}

\newcommand{\defi}[1]{\textsf{#1}} 

\newcommand{\C}{\mathbb{C}}
\newcommand{\F}{\mathbb{F}}

\newcommand{\ii}{\mathbf{i}}

\newcommand{\Q}{\mathbb{Q}}
\newcommand{\R}{\mathbb{R}}
\newcommand{\xx}{\mathbf{x}}
\newcommand{\Z}{\mathbb{Z}}

\newcommand{\mm}{\mathfrak{m}}

\newcommand{\tH}{{\operatorname{th}}}

\newcommand{\GL}{\operatorname{GL}}

\newtheorem{theorem}{Theorem}

\theoremstyle{definition}

\theoremstyle{remark}
\newtheorem{remark}[theorem]{Remark}

\usepackage[
	pdfauthor={Bjorn Poonen}, 
]{hyperref}
\usepackage[alphabetic,lite]{amsrefs} 

\begin{document}

\title{$p$-adic interpolation of iterates}
\author{Bjorn Poonen}
\thanks{This research was supported by National Science Foundation grant DMS-1069236.}
\address{Department of Mathematics, Massachusetts Institute of Technology, Cambridge, MA 02139-4307, USA}
\email{poonen@math.mit.edu}
\urladdr{\url{http://math.mit.edu/~poonen/}}
\date{July 22, 2013}

\begin{abstract}
Extending work of Bell and of Bell, Ghioca, and Tucker,
we prove that for a $p$-adic analytic self-map $f$ on a
closed unit polydisk,
if every coefficient of $f(\xx)-\xx$ has valuation
greater than that of $p^{1/(p-1)}$, 
then the iterates of $f$ can be $p$-adically interpolated;
i.e., there exists a function $g(\xx,n)$ analytic in both $\xx$ and $n$
such that $g(\xx,n) = f^n(\xx)$ whenever $n \in \Z_{\ge 0}$.
\end{abstract}

\maketitle


Inspired by the work of Skolem~\cite{Skolem1934}, 
Mahler~\cites{Mahler1935},
and Lech~\cite{Lech1953}
on linear recursive sequences,
Bell~\cite{Bell2008} 
proved that for a suitable $p$-adic analytic function $f$
and starting point $\xx$, the iterate-computing map $n \mapsto f^n(\xx)$ 
extends to a $p$-adic analytic function $g(n)$ defined for $n \in \Z_p$.
This result, along with its generalization by Bell, Ghioca, and Tucker 
in~\cite{Bell-Ghioca-Tucker2010}*{\S3}
and earlier linearization results 
by Herman and Yoccoz \cite{Herman-Yoccoz1983}*{Theorem~1} 
and Rivera-Letelier \cite{Rivera-Letelier2003}*{\S 3.2},
has significance beyond its intrinsic interest,
because of its applications towards
the dynamical Mordell--Lang conjecture
\citelist{
\cite{Bell2006} 
\cite{Ghioca-Tucker2009} 
\cite{Bell-Ghioca-Tucker2010} 
\cite{Benedetto-Ghioca-Kurlberg-Tucker2012} 
\cite{Benedetto-Ghioca-Hutz-Kurlberg-Scanlon-Tucker2013}
}.

Our main result, Theorem~\ref{T:main},
is a variant that is best possible
(in a sense explained in Remark~\ref{R:best possible}).
Our proof is new even over $\Q_p$, 
and extends immediately to more general valued fields.
It settles an open question about the case $p=3$.
The function $g$ we obtain is analytic in $\xx$ as well as $n$.

We now set the notation for our statement.
Let $p$ be a prime number.
Let $K$ be a field 
that is complete with respect to an absolute value $|\;|$ satisfying $|p|=1/p$.
Let $R$ be the valuation ring in $K$.
For $f \in R[\xx] \colonequals R[x_1,\ldots,x_d]$, 
let $\|f\|$ 
be the supremum of the absolute values of the coefficients of $f$.
The \defi{Tate algebra} $R\langle \xx \rangle$ is the completion of $R[\xx]$
with respect to $\|\;\|$.
More concretely, $R\langle \xx \rangle$ 
is the set of $f = \sum_{\ii \in \Z_{\ge 0}^d} f_{\ii} \xx^{\ii} \in R[[\xx]]$
converging on the closed unit polydisk;
convergence is equivalent to $|f_{\ii}| \to 0$ as $\ii \to \infty$.
For $f,g \in R\langle \xx \rangle$ and $c \in \R_{\ge 0}$,
the notation $f \in p^c R\langle \xx \rangle$ means $\|f\| \le |p|^c$,
and $f \equiv g \pmod{p^c}$ means $\|f-g\| \le |p|^c$;
extend componentwise to $f,g \in R\langle \xx \rangle^d$.

\begin{theorem}
\label{T:main}
If $f \in R\langle x_1,\ldots,x_d \rangle^d$
satisfies $f(\xx) \equiv \xx \pmod{p^c}$ for some $c>\frac{1}{p-1}$,
then there exists $g \in R\langle x_1,\ldots,x_d,n \rangle$
such that $g(\xx,n) = f^n(\xx)$ in $R\langle \xx \rangle^d$ 
for each $n \in \Z_{\ge 0}$.
\end{theorem}

Our proof will check directly that the Mahler series~\cite{Mahler1958} 
interpolating the sequence
\[
	\xx, f(\xx), f(f(\xx)), \ldots
\]
converges to an analytic function.
This is the difference operator analogue 
of proving that a function $\phi$ is analytic
by checking that its Taylor series converges to $\phi$.

\begin{proof}
Since $f(\xx) \equiv \xx \pmod{p^c}$,
we have $h(f(\xx)) \equiv h(\xx) \pmod{p^c}$
for any $h \in R[\xx]^d$ 
and (by taking limits) also for any $h \in R\langle \xx \rangle^d$.
In other words, the linear operator $\Delta$ defined by
\[
	(\Delta h)(\xx) \colonequals h(f(\xx)) - h(\xx)
\]
maps $R\langle \xx \rangle^d$ into $p^c R\langle \xx \rangle^d$.
In particular, $m$ applications of $\Delta$ to the identity function 
yields $\Delta^m \xx \in p^{mc} R\langle \xx \rangle^d$.
On the other hand, $|m!| \ge p^{-m/(p-1)}$. 
Thus the Mahler series
\[
	g(\xx,n) \colonequals \sum_{m \ge 0} \binom{n}{m} \Delta^m \xx
	= \sum_{m \ge 0} n(n-1)\cdots(n-m+1) \frac{ \Delta^m \xx}{m!}
\]
converges in $R\langle \xx,n \rangle^d$ with respect to $\|\;\|$.
Let $I$ be the identity operator.
If $n \in \Z_{\ge 0}$, then 
\[
	g(\xx,n) = \sum_{m=0}^n \binom{n}{m} \Delta^m \xx 
	= (\Delta + I)^n \xx = f^n(\xx).\qedhere
\]
\end{proof}

\begin{remark}
The relation $g(\xx,n+1)=f(g(\xx,n))$ in $R\langle \xx \rangle^d$ 
holds for each $n$ in the infinite set $\Z_{\ge 0}$,
so it is an identity in $R \langle \xx,n \rangle^d$.
\end{remark}

\begin{remark}
\label{R:best possible}
The hypothesis on $f$ holds for $K=\Q_p$
if $f(\xx) \equiv \xx \pmod{p}$ and $p \ge 3$;
previously the conclusion was known only for $p \ge 5$ 
\citelist{\cite{Bell2008} \cite{Bell-Ghioca-Tucker2010}*{\S3}}.
On the other hand, $f(x) \colonequals -x$ is a counterexample for $p=2$ 
\cite{Bell2008}*{\S3}.
Similarly, the inequality on $c$ in Theorem~\ref{T:main} is best possible
for each $p$: consider $f(x) \colonequals \zeta x$ where $\zeta$
is a primitive $p^{\tH}$ root of unity in $\C_p$.
\end{remark}

\begin{remark}
Let $\mm$ be the maximal ideal of $R$.
Let $k \colonequals R/\mm$.
If $f(\xx) \bmod \mm = \xx$,
so that $f(\xx) \equiv \xx \pmod{p^c}$ holds for some $c>0$,
then $f^p(\xx) \equiv \xx \pmod{p^c}$ holds for a larger $c$,
and by iterating we find $r \in \Z_{\ge 0}$ 
such that Theorem~\ref{T:main} applies to $f^{p^r}$.
More generally, 
if $f(\xx) \bmod \mm = A\xx$ for some $A \in \GL_d(k)$ of finite order,
then there exists $s \in \Z_{>0}$ such that 
$f^s$ satisfies the hypothesis of Theorem~\ref{T:main}.
This finite order hypothesis is automatic if $K$ is $\Q_p$ or $\C_p$
since then $k$ is algebraic over $\F_p$
and \emph{every} element of $\GL_d(k)$ is of finite order.
Cf.~\cite{Bell-Ghioca-Tucker2010}*{\S2.2}.
\end{remark}


\section*{Acknowledgement} 

We thank Thomas J. Tucker for a helpful suggestion on the exposition.

\end{document}